\documentclass[11pt,article]{amsart}
\UseRawInputEncoding

\usepackage{amsmath,amssymb,color,accents}
\usepackage{graphicx,color}
\usepackage[colorlinks=true, pdfstartview=FitV, linkcolor=magenta, citecolor=magenta, urlcolor=blue]{hyperref}
\usepackage[nameinlink,capitalise]{cleveref}
\usepackage{comment}
\usepackage{mathtools}
\usepackage[all]{xy}
\usepackage{ytableau}
\usepackage{amsfonts}
\usepackage{longtable}
\usepackage{enumerate}
\usepackage{tikz}
\usepackage{multicol}
\usepackage{tabularx}
\usepackage[justification=centering]{caption}

\usepackage{fancyvrb}
\usepackage{fvextra}
\usepackage{listings}
\usepackage{listings}
\usepackage{xcolor}
\definecolor{maccolor}{rgb}{0.3,0.3,0.8}
\lstdefinelanguage{Macaulay2}
{
basicstyle={\ttfamily},
keywordstyle={\color{maccolor!80!black}},
commentstyle={\color{gray}},
stringstyle={\color{red!40!black}},
rulecolor=\color{maccolor},
basewidth={1.2ex}, 
sensitive=false,
morestring=[b]",
escapechar={`},
escapebegin={\rmfamily},
morekeywords={
ConnectionMatrices, connectionMatrix, gaugeMatrix, gaugeTransform, holonomicRank, isEpsilonFactorized, isIntegrable, makeWeylAlgebra, standardMonomials, baseFractionField,
about,abs,AbstractToricVarieties,accumulate,Acknowledgement,acos,acosh,acot,addCancelTask,addDependencyTask,addEndFunction,addHook,AdditionalPaths,addStartFunction,addStartTask,Adjacent,adjoint,AdjointIdeal,AffineVariety,AfterEval,AfterNoPrint,AfterPrint,agm,AInfinity,alarm,AlgebraicSplines,Algorithm,Alignment,AllCodimensions,allowableThreads,ambient,analyticSpread,Analyzer,AnalyzeSheafOnP1,ancestor,ancestors,ANCHOR,and,andP,AngleBarList,ann,annihilator,antipode,any,append,applicationDirectory,applicationDirectorySuffix,apply,applyKeys,applyPairs,applyTable,applyValues,apropos,argument,Array,arXiv,Ascending,ascii,asin,asinh,ass,assert,associatedGradedRing,associatedPrimes,AssociativeAlgebras,AssociativeExpression,atan,atan2,atEndOfFile,Authors,autoload,AuxiliaryFiles,backtrace,Bag,Bareiss,baseFilename,BaseFunction,baseName,baseRing,baseRings,BaseRow,BasicList,basis,BasisElementLimit,Bayer,BeforePrint,beginDocumentation,BeginningMacaulay2,Benchmark,benchmark,Bertini,BesselJ,BesselY,betti,BettiCharacters,BettiTally,between,BGG,BIBasis,Binary,BinaryOperation,Binomial,binomial,BinomialEdgeIdeals,Binomials,BKZ,BlockMatrix,BLOCKQUOTE,BODY,Body,BoijSoederberg,BOLD,Book3264Examples,Boolean,BooleanGB,borel,Boxes,BR,break,Browse,Bruns,cache,CacheExampleOutput,CacheFunction,CacheTable,cacheValue,CallLimit,cancelTask,capture,catch,Caveat,CC,CDATA,ceiling,Center,centerString,Certification,ChainComplex,chainComplex,ChainComplexExtras,ChainComplexMap,ChainComplexOperations,ChangeMatrix,char,CharacteristicClasses,characters,charAnalyzer,check,CheckDocumentation,chi,Chordal,class,Classic,clean,clearAll,clearEcho,clearOutput,close,closeIn,closeOut,ClosestFit,CODE,code,codim,CodimensionLimit,coefficient,CoefficientRing,coefficientRing,Cofactor,CohenEngine,CohenTopLevel,CoherentSheaf,CohomCalg,cohomology,coimage,CoincidentRootLoci,coker,cokernel,collectGarbage,columnAdd,columnate,columnMult,columnPermute,columnRankProfile,columnSwap,combine,Command,commandInterpreter,commandLine,COMMENT,commonest,commonRing,comodule,CompactMatrix,compactMatrixForm,CompiledFunction,CompiledFunctionBody,CompiledFunctionClosure,Complement,complement,complete,CompleteIntersection,CompleteIntersectionResolutions,Complexes,ComplexField,components,compose,compositions,compress,concatenate,conductor,ConductorElement,cone,ConformalBlocks,conjugate,connectionCount,Consequences,Constant,Constants,constParser,content,continue,contract,Contributors,ConvexInterface,conwayPolynomial,ConwayPolynomials,copy,copyDirectory,copyFile,copyright,Core,CorrespondenceScrolls,cos,cosh,cot,CotangentSchubert,cotangentSheaf,coth,cover,coverMap,cpuTime,createTask,Cremona,csc,csch,current,currentColumnNumber,currentDirectory,currentFileDirectory,currentFileName,currentLayout,currentLineNumber,currentPackage,currentString,currentTime,Cyclotomic,Database,Date,DD,dd,deadParser,debug,debugError,DebuggingMode,debuggingMode,debugLevel,DecomposableSparseSystems,Decompose,decompose,deepSplice,Default,default,defaultPrecision,Degree,degree,degreeLength,DegreeLift,DegreeLimit,DegreeMap,DegreeOrder,DegreeRank,Degrees,degrees,degreesMonoid,degreesRing,delete,demark,denominator,Dense,Density,Depth,depth,Descending,Descent,Describe,describe,Description,det,determinant,DeterminantalRepresentations,DGAlgebras,diagonalMatrix,diameter,Dictionary,dictionary,dictionaryPath,diff,DiffAlg,difference,dim,directSum,disassemble,discriminant,dismiss,Dispatch,distinguished,DIV,Divide,divideByVariable,DivideConquer,DividedPowers,Divisor,DL,Dmodules,do,doc,docExample,docTemplate,document,DocumentTag,Down,drop,DT,dual,eagonNorthcott,EagonResolution,echoOff,echoOn,EdgeIdeals,edit,EigenSolver,eigenvalues,eigenvectors,eint,EisenbudHunekeVasconcelos,elapsedTime,elapsedTiming,elements,Eliminate,eliminate,Elimination,EliminationMatrices,EllipticCurves,EllipticIntegrals,else,EM,Email,End,end,endl,endPackage,Engine,engineDebugLevel,EngineRing,EngineTests,entries,EnumerationCurves,environment,Equation,EquivariantGB,erase,erf,erfc,error,errorDepth,euler,EulerConstant,eulers,even,EXAMPLE,ExampleFiles,ExampleItem,examples,ExampleSystems,Exclude,exec,exit,exp,expectedReesIdeal,expm1,exponents,export,exportFrom,exportMutable,Expression,expression,Ext,extend,ExteriorIdeals,ExteriorModules,exteriorPower,Factor,factor,false,Fano,FastMinors,FastNonminimal,FGLM,File,fileDictionaries,fileExecutable,fileExists,fileExitHooks,fileLength,fileMode,FileName,FilePosition,fileReadable,fileTime,fileWritable,fillMatrix,findFiles,findHeft,FindOne,findProgram,findSynonyms,FiniteFittingIdeals,First,first,firstkey,FirstPackage,fittingIdeal,flagLookup,FlatMonoid,flatten,flattenRing,Flexible,flip,floor,flush,fold,FollowLinks,for,forceGB,fork,FormalGroupLaws,Format,format,formation,FourierMotzkin,FourTiTwo,fpLLL,frac,fraction,FractionField,frames,FrobeniusThresholds,from,fromDividedPowers,fromDual,FunctionApplication,FunctionBody,functionBody,FunctionClosure,FunctionFieldDesingularization,fusePairs,futureParser,GaloisField,Gamma,gb,GBDegrees,gbRemove,gbSnapshot,gbTrace,gcd,gcdCoefficients,gcdLLL,GCstats,genera,GeneralOrderedMonoid,GenerateAssertions,generateAssertions,generator,generators,Generic,GenericInitialIdeal,genericMatrix,genericSkewMatrix,genericSymmetricMatrix,gens,genus,getc,getChangeMatrix,getenv,getGlobalSymbol,getNetFile,getNonUnit,getPrimeWithRootOfUnity,getSymbol,getWWW,GF,gfanInterface,Givens,GKMVarieties,GLex,Global,global,globalAssign,globalAssignFunction,GlobalAssignHook,globalAssignment,globalAssignmentHooks,GlobalDictionary,GlobalHookStore,globalReleaseFunction,GlobalReleaseHook,Gorenstein,GradedLieAlgebras,GradedModule,gradedModule,GradedModuleMap,gradedModuleMap,gramm,GraphicalModels,GraphicalModelsMLE,Graphics,graphIdeal,graphRing,Graphs,Grassmannian,GRevLex,GroebnerBasis,groebnerBasis,GroebnerBasisOptions,GroebnerStrata,GroebnerWalk,groupID,GroupLex,GroupRevLex,GTZ,Hadamard,handleInterrupts,HardDegreeLimit,hash,HashTable,hashTable,HEAD,HEADER1,HEADER2,HEADER3,HEADER4,HEADER5,HEADER6,HeaderType,Heading,Headline,Heft,heft,Height,height,help,Hermite,hermite,Hermitian,HH,hh,HigherCIOperators,HighestWeights,Hilbert,hilbertFunction,hilbertPolynomial,hilbertSeries,HodgeIntegrals,hold,Holder,Hom,homeDirectory,HomePage,Homogeneous,Homogeneous2,homogenize,homology,homomorphism,HomotopyLieAlgebra,hooks,horizontalJoin,HorizontalSpace,HR,HREF,HTML,html,httpHeaders,Hybrid,HyperplaneArrangements,Hypertext,hypertext,HypertextContainer,HypertextParagraph,icFracP,icFractions,icMap,icPIdeal,id,Ideal,ideal,idealizer,identity,if,IgnoreExampleErrors,ii,image,imaginaryPart,IMG,ImmutableType,importFrom,incomparable,Increment,independentSets,indeterminate,IndeterminateNumber,Index,index,indexComponents,IndexedVariable,IndexedVariableTable,indices,inducedMap,inducesWellDefinedMap,InexactField,InexactFieldFamily,InexactNumber,InfiniteNumber,infinity,info,InfoDirSection,infoHelp,Inhomogeneous,input,Inputs,insert,installAssignmentMethod,installedPackages,installHilbertFunction,installMethod,installMinprimes,installPackage,InstallPrefix,instance,instances,IntegralClosure,integralClosure,integrate,IntermediateMarkUpType,interpreterDepth,intersect,intersectInP,interval,InvariantRing,inverse,InverseMethod,inversePermutation,Inverses,inverseSystem,InverseSystems,Invertible,InvolutiveBases,irreducibleCharacteristicSeries,irreducibleDecomposition,isAffineRing,isANumber,isBorel,isCanceled,isCommutative,isConstant,isDirectory,isDirectSum,isEmpty,isField,isFinite,isFinitePrimeField,isFreeModule,isGlobalSymbol,isHomogeneous,isIdeal,isInfinite,isInjective,isInputFile,isIsomorphism,isLinearType,isListener,isLLL,isMember,isModule,isMonomialIdeal,isNormal,isOpen,isOutputFile,isPolynomialRing,isPrimary,isPrime,isPrimitive,isPseudoprime,isQuotientModule,isQuotientOf,isQuotientRing,isReady,isReal,isReduction,isRegularFile,isRing,isSkewCommutative,isSorted,isSquareFree,isStandardGradedPolynomialRing,isSubmodule,isSubquotient,isSubset,isSupportedInZeroLocus,isSurjective,isTable,isUnit,isWellDefined,isWeylAlgebra,ITALIC,Iterate,Jacobian,jacobian,jacobianDual,Jets,Join,join,Jupyter,K3Carpets,K3Surfaces,Keep,KeepFiles,KeepZeroes,ker,kernel,kernelLLL,kernelOfLocalization,Key,keys,Keyword,Keywords,kill,koszul,Kronecker,KustinMiller,LABEL,last,lastMatch,LATER,LatticePolytopes,Layout,lcm,leadCoefficient,leadComponent,leadMonomial,leadTerm,Left,left,length,LengthLimit,letterParser,Lex,LexIdeals,LI,Licenses,LieTypes,lift,liftable,Limit,limitFiles,limitProcesses,LinearAlgebra,LinearTruncations,lineNumber,LINK,linkFile,List,list,listForm,listLocalSymbols,listSymbols,listUserSymbols,LITERAL,LLL,LLLBases,lngamma,load,loadDepth,LoadDocumentation,loadedFiles,loadedPackages,loadPackage,Local,local,localDictionaries,LocalDictionary,localize,LocalRings,locate,log,log1p,LongPolynomial,lookup,lookupCount,LowerBound,LUdecomposition,M0nbar,M2CODE,Macaulay2Doc,makeDirectory,MakeDocumentation,makeDocumentTag,MakeHTML,MakeInfo,MakeLinks,makePackageIndex,MakePDF,makeS2,Manipulator,map,MapExpression,MapleInterface,markedGB,Markov,MarkUpType,match,mathML,Matrix,matrix,MatrixExpression,Matroids,max,maxAllowableThreads,maxExponent,MaximalRank,maxPosition,MaxReductionCount,MCMApproximations,member,memoize,memoizeClear,memoizeValues,MENU,merge,mergePairs,META,method,MethodFunction,MethodFunctionBinary,MethodFunctionSingle,MethodFunctionWithOptions,methodOptions,methods,midpoint,min,minExponent,mingens,mingle,minimalBetti,MinimalGenerators,MinimalMatrix,minimalPresentation,minimalPresentationMap,minimalPresentationMapInv,MinimalPrimes,minimalPrimes,minimalReduction,Minimize,minimizeFilename,MinimumVersion,minors,minPosition,minPres,minprimes,Minus,minus,Miura,MixedMultiplicity,mkdir,mod,Module,module,ModuleDeformations,modulo,MonodromySolver,Monoid,monoid,MonoidElement,Monomial,MonomialAlgebras,monomialCurveIdeal,MonomialIdeal,monomialIdeal,MonomialIntegerPrograms,MonomialOrbits,MonomialOrder,Monomials,monomials,MonomialSize,monomialSubideal,moveFile,multidegree,multidoc,multigraded,MultigradedBettiTally,MultiGradedRationalMap,MultiplicitySequence,MultiplierIdeals,MultiplierIdealsDim2,MultiprojectiveVarieties,mutable,MutableHashTable,mutableIdentity,MutableList,MutableMatrix,mutableMatrix,NAGtypes,Name,nanosleep,Nauty,NautyGraphs,NCAlgebra,NCLex,needs,needsPackage,Net,net,NetFile,netList,new,newClass,newCoordinateSystem,NewFromMethod,newline,NewMethod,newNetFile,NewOfFromMethod,NewOfMethod,newPackage,newRing,nextkey,nextPrime,nil,NNParser,NoetherianOperators,NoetherNormalization,NonminimalComplexes,nonspaceAnalyzer,NoPrint,norm,normalCone,Normaliz,NormalToricVarieties,not,Nothing,notify,notImplemented,NTL,null,nullaryMethods,nullhomotopy,nullParser,nullSpace,Number,number,NumberedVerticalList,numcols,numColumns,numerator,numeric,NumericalAlgebraicGeometry,NumericalCertification,NumericalImplicitization,NumericalLinearAlgebra,NumericalSchubertCalculus,numericInterval,NumericSolutions,numgens,numRows,numrows,odd,oeis,ofClass,OL,OldPolyhedra,OldToricVectorBundles,on,OneExpression,OnlineLookup,OO,oo,ooo,oooo,openDatabase,openDatabaseOut,openFiles,openIn,openInOut,openListener,OpenMath,openOut,openOutAppend,operatorAttributes,Option,OptionalComponentsPresent,optionalSignParser,Options,options,OptionTable,optP,or,Order,order,OrderedMonoid,orP,OutputDictionary,Outputs,override,pack,Package,package,PackageCitations,PackageDictionary,PackageExports,PackageImports,PackageTemplate,packageTemplate,pad,pager,PairLimit,pairs,PairsRemaining,PARA,Parametrization,parent,Parenthesize,Parser,Parsing,part,Partition,partition,partitions,parts,path,pdim,peek,PencilsOfQuadrics,Permanents,permanents,permutations,pfaffians,PHCpack,PhylogeneticTrees,pi,PieriMaps,pivots,PlaneCurveSingularities,plus,poincare,poincareN,polarize,poly,Polyhedra,Polymake,PolynomialRing,Posets,Position,position,positions,PositivityToricBundles,POSIX,Postfix,Power,power,powermod,PRE,Precision,precision,Prefix,prefixDirectory,prefixPath,preimage,prepend,presentation,pretty,primaryComponent,PrimaryDecomposition,primaryDecomposition,PrimaryTag,PrimitiveElement,Print,print,printerr,printingAccuracy,printingLeadLimit,printingPrecision,printingSeparator,printingTimeLimit,printingTrailLimit,printString,printWidth,processID,Product,product,ProductOrder,profile,profileSummary,Program,programPaths,ProgramRun,Proj,Projective,ProjectiveHilbertPolynomial,projectiveHilbertPolynomial,ProjectiveVariety,promote,protect,Prune,prune,PruneComplex,pruningMap,Pseudocode,pseudocode,pseudoRemainder,Pullback,PushForward,pushForward,Python,QQ,QQParser,QRDecomposition,QthPower,Quasidegrees,QuaternaryQuartics,QuillenSuslin,quit,Quotient,quotient,quotientRemainder,QuotientRing,Radical,radical,RadicalCodim1,radicalContainment,RaiseError,random,RandomCanonicalCurves,RandomComplexes,RandomCurves,RandomCurvesOverVerySmallFiniteFields,RandomGenus14Curves,RandomIdeals,randomKRationalPoint,RandomMonomialIdeals,randomMutableMatrix,RandomObjects,RandomPlaneCurves,RandomPoints,RandomSpaceCurves,Range,rank,RationalMaps,RationalPoints,RationalPoints2,ReactionNetworks,read,readDirectory,readlink,readPackage,RealField,RealFP,realPart,realpath,RealQP,RealQP1,RealRoots,RealRR,RealXD,recursionDepth,recursionLimit,Reduce,reducedRowEchelonForm,reduceHilbert,reductionNumber,ReesAlgebra,reesAlgebra,reesAlgebraIdeal,reesIdeal,References,ReflexivePolytopesDB,regex,regexQuote,registerFinalizer,regSeqInIdeal,Regularity,regularity,relations,RelativeCanonicalResolution,relativizeFilename,Reload,remainder,RemakeAllDocumentation,removeDirectory,removeFile,removeLowestDimension,reorganize,replace,RerunExamples,res,reshape,ResidualIntersections,ResLengthThree,Resolution,resolution,ResolutionsOfStanleyReisnerRings,restart,Result,resultant,Resultants,return,returnCode,Reverse,reverse,RevLex,Right,right,Ring,ring,RingElement,RingFamily,ringFromFractions,RingMap,rootPath,roots,rootURI,rotate,round,rowAdd,RowExpression,rowMult,rowPermute,rowRankProfile,rowSwap,RR,RRi,rsort,run,RunDirectory,RunExamples,RunExternalM2,runHooks,runLengthEncode,runProgram,same,saturate,Saturation,scan,scanKeys,scanLines,scanPairs,scanValues,schedule,schreyerOrder,Schubert,Schubert2,SchurComplexes,SchurFunctors,SchurRings,SCRIPT,scriptCommandLine,ScriptedFunctor,SCSCP,searchPath,sec,sech,SectionRing,SeeAlso,seeParsing,SegreClasses,select,selectInSubring,selectVariables,SelfInitializingType,SemidefiniteProgramming,Seminormalization,separate,SeparateExec,separateRegexp,Sequence,sequence,Serialization,serialNumber,Set,set,setEcho,setGroupID,setIOExclusive,setIOSynchronized,setIOUnSynchronized,setRandomSeed,setup,setupEmacs,sheaf,SheafExpression,sheafExt,sheafHom,SheafOfRings,shield,ShimoyamaYokoyama,short,show,showClassStructure,showHtml,showStructure,showTex,showUserStructure,SimpleDoc,simpleDocFrob,SimplicialComplexes,SimplicialDecomposability,SimplicialPosets,SimplifyFractions,sin,singularLocus,sinh,size,size2,SizeLimit,SkewCommutative,SlackIdeals,sleep,SLnEquivariantMatrices,SLPexpressions,SMALL,smithNormalForm,solve,someTerms,Sort,sort,sortColumns,SortStrategy,source,SourceCode,SourceRing,SPACE,SpaceCurves,SPAN,span,SparseMonomialVectorExpression,SparseResultants,SparseVectorExpression,Spec,SpechtModule,SpecialFanoFourfolds,specialFiber,specialFiberIdeal,SpectralSequences,splice,splitWWW,sqrt,SRdeformations,stack,stacksProject,Standard,standardForm,standardPairs,StartWithOneMinor,stashValue,StatePolytope,StatGraphs,status,stderr,stdio,step,StopBeforeComputation,stopIfError,StopWithMinimalGenerators,Strategy,String,STRONG,StronglyStableIdeals,STYLE,Style,style,SUB,sub,SubalgebraBases,sublists,submatrix,submatrixByDegrees,Subnodes,subquotient,SubringLimit,Subscript,subscript,SUBSECTION,subsets,substitute,substring,subtable,Sugarless,Sum,sum,SumOfTwists,SumsOfSquares,SUP,super,SuperLinearAlgebra,Superscript,superscript,support,SVD,SVDComplexes,switch,SwitchingFields,sylvesterMatrix,Symbol,symbol,SymbolBody,symbolBody,SymbolicPowers,symlinkDirectory,symlinkFile,symmetricAlgebra,symmetricAlgebraIdeal,symmetricKernel,SymmetricPolynomials,symmetricPower,synonym,SYNOPSIS,syz,Syzygies,SyzygyLimit,SyzygyMatrix,SyzygyRows,syzygyScheme,TABLE,Table,table,take,Tally,tally,tan,TangentCone,tangentCone,tangentSheaf,tanh,target,Task,taskResult,TateOnProducts,TD,temporaryFileName,tensor,tensorAssociativity,TensorComplexes,terminalParser,terms,TEST,Test,testExample,testHunekeQuestion,TestIdeals,TestInput,tests,TEX,tex,TeXmacs,texMath,Text,TH,then,Thing,ThinSincereQuivers,ThreadedGB,threadVariable,Threshold,throw,Time,time,times,timing,TITLE,TO,to,TO2,toAbsolutePath,toCC,toDividedPowers,toDual,toExternalString,toField,TOH,toList,toLower,top,top,topCoefficients,Topcom,topComponents,topLevelMode,Tor,TorAlgebra,Toric,ToricInvariants,ToricTopology,ToricVectorBundles,toRR,toRRi,toSequence,toString,TotalPairs,toUpper,TR,trace,transpose,TriangularSets,Tries,Trim,trim,Triplets,Tropical,true,Truncate,truncate,truncateOutput,Truncations,try,TSpreadIdeals,TT,tutorial,Type,TypicalValue,typicalValues,UL,ultimate,unbag,uncurry,Undo,undocumented,uniform,uninstallAllPackages,uninstallPackage,Unique,unique,Units,Unmixed,unsequence,unstack,Up,UpdateOnly,UpperTriangular,URL,urlEncode,Usage,use,UseCachedExampleOutput,UseHilbertFunction,UserMode,userSymbols,UseSyzygies,utf8,utf8check,validate,value,values,Variable,VariableBaseName,Variables,Variety,variety,vars,Vasconcelos,Vector,vector,VectorExpression,VectorFields,VectorGraphics,Verbose,Verbosity,Verify,VersalDeformations,versalEmbedding,Version,version,VerticalList,VerticalSpace,viewHelp,VirtualResolutions,VirtualTally,VisibleList,Visualize,wait,WebApp,wedgeProduct,weightRange,Weights,WeylAlgebra,WeylGroups,when,whichGm,while,width,wikipedia,Wrap,wrap,WrapperType,XML,xor,youngest,zero,ZeroExpression,zeta,ZZ,ZZParser}
}
\lstalias{Macaulay2output}{Macaulay2}

\usepackage[inline]{enumitem} 

\usepackage{tikz}
\usepackage[all]{xy}
\usetikzlibrary{matrix,cd}

\tikzset{
  treenode/.style = {align=center, inner sep=0pt, text centered,solid,thin,
    font=\sffamily},
  arn_n/.style = {treenode, circle, white, font=\sffamily\bfseries, draw=black,
    fill=black, text width=.5em},
  arn_nl/.style = {treenode, circle, white, font=\sffamily\bfseries, draw=black,
    fill=black, text width=1.5em},  
  arn_r/.style = {treenode, circle, red, draw=red, 
    text width=.5em, very thick},
  arn_v/.style = {treenode, circle, black, font=\sffamily\bfseries, draw=black, text width=1.2em},
  arn_x/.style = {treenode, rectangle, draw=black,
    minimum width=.5em, minimum height=0.5em},
  dott/.style={edge from parent/.style={dotted, very thick,circle,draw}},
  emph/.style={edge from parent/.style={dashed, very thick,circle,draw}},
  norm/.style={edge from parent/.style={solid,thin,circle,draw}}
}

\usepackage{tabstackengine}

\usepackage{listings}

\let\amsamp=&
\usepackage{bbm}

\newcommand{\sbm}[1]{{\let\amp=&\begin{pmatrix}#1\end{pmatrix}}}

\numberwithin{equation}{section}
\newtheorem{theorem}{Theorem}[section]
\newtheorem{lemma}[theorem]{Lemma}
\newtheorem{proposition}[theorem]{Proposition}
\newtheorem{corollary}[theorem]{Corollary}
\newtheorem{conjecture}[theorem]{Conjecture}

\newtheorem*{theorem*}{Theorem}

\theoremstyle{definition}
\newtheorem{definition}[theorem]{Definition}

\theoremstyle{remark}
\newtheorem{remark}[theorem]{Remark}
\newtheorem{example}[theorem]{Example}

\newtheorem{observation}[theorem]{Observation}
\newtheorem{acknowledgments}[theorem]{Acknowledgments}

\newcommand{\bpf}{\begin{proof}}
\newcommand{\epf}{\end{proof}}
\newcommand{\bpr}{\begin{proposition}}
\newcommand{\epr}{\end{proposition}}
\newcommand{\bdf}{\begin{definition}}
\newcommand{\edf}{\end{definition}}
\newcommand{\blm}{\begin{lemma}}
\newcommand{\elm}{\end{lemma}}
\newcommand{\bex}{\begin{example}\rm }
\newcommand{\eex}{\end{example}}
\newcommand{\bcor}{\begin{corollary}}
\newcommand{\ecor}{\end{corollary}}
\newcommand{\bthm}{\begin{theorem}}
\newcommand{\ethm}{\end{theorem}}
\newcommand{\be}{\begin{enumerate}}
\newcommand{\ee}{\end{enumerate}}
\newcommand{\bq}{\begin{equation}}
\newcommand{\eq}{\end{equation}}
\newcommand{\bb}{\begin{itemize}}
\newcommand{\eb}{\end{itemize}}
\newcommand{\bpw}{\begin{cases}}
\newcommand{\epw}{\end{cases}}

\DeclareMathOperator{\ev}{ev}
\DeclareMathOperator{\Proj}{Proj}

\DeclareMathOperator{\Span}{span}

\newcommand{\cB}{{\mathcal B}}

\newcommand{\cL}{{\mathcal L}}

\newcommand{\cO}{{\mathcal O}}

\newcommand{\C}{{\mathbb C}}

\newcommand{\Q}{{\mathbb Q}}
\newcommand{\R}{{\mathbb R}}

\newcommand{\Z}{{\mathbb Z}}
\renewcommand{\AA}{{\mathbb A}}
\newcommand{\PP}{{\mathbb P}}
\numberwithin{equation}{section}

\newcommand{\LF}{\left\lfloor}
\newcommand{\RF}{\right\rfloor}

\title{On the algebraic properties of the B\"or\"oczky configuration}

\author[Jake Kettinger]{Jake Kettinger}
\address{Department of Mathematics, Colorado State University, 841 Oval Drive, Fort Collins, CO, 80521}
\email{jake.kettinger@colostate.edu}

\author{Shahriyar Roshan Zamir}
\address{Department of Mathematics, Tulane University, 6823 St. Charles Avenue,
New Orleans, LA, 70118}
\email{sroshanzamir@tulane.edu}

\usepackage{fancyhdr}
\begin{document}

\begin{abstract}
The B\"or\"oczky configuration of lines and (multiple) points exhibits extremal behavior in commutative algebra and combinatorics. Examples of this appear in the context of the containment problem for ordinary and symbolic powers and the proof of the Dirac-Motzkin conjecture by Green and Tao. This paper studies the algebraic properties of B\"or\"oczky configurations for arbitrary values of $n$. 
Our results compute the Waldschmit constant of the defining ideal of these configurations. Moreover, we use the weighted projective plane $\PP(1,2,3)$ to give an upper bound for the degree of the minimal generators of their ideal. Finally, this construction is applied to an elliptic curve in $\mathbb{P}^2$ to give a new counterexample to the containment $I^{(3)}\subseteq I^2$.
\end{abstract}

\maketitle

\setcounter{tocdepth}{1} 
\tableofcontents

\section{Introduction}
The goal of this note is to study the algebraic properties of the B\"or\"oczky configuration of lines and points. The B\"or\"oczky configuration, defined for every $n\geq 3$, consists of lines, ordinary points, where exactly two lines meet, and triple points, where exactly three lines of the configuration meet; see \Cref{sec: Borocky General} for a precise definition. Let $B_n$ denote the set of all lines and all intersection points of the B\"or\"oczky configuration. 

There are two competing definitions of this configuration in the literature. One used by Green and Tao \cite{greentao} in their proof of the Dirac-Motzkin conjecture, and one used by F\"uredi and Pal\'asti \cite{PalastiFerudi} which, among other results, provides a negative answer to a question of Paul Erd\"os. An essential difference between the two is that Green and Tao consider this construction in the real {\it projective} plane while F\"uredi and Pal\'asti make use of it in the real {\it affine} plane. The latter construction described in \cite{PalastiFerudi}, considered in the affine chart $U_2$ of $\mathbb{P}^2$, is the one of interest to us and henceforth is referred to ``the B\"or\"oczky configuration.'' It it notable that their construction is a special case of an earlier construction by Burr, Gr\"unbaum and Sloane \cite{GrunbaumSloane} which uses elliptic integration; see \cite[Section 6]{PalastiFerudi}. We draw another connection to elliptic curves in \Cref{sec: elliptified Boroczky}.

The algebraic context in which the B\"or\"oczky configuration arises is the containment problem of ordinary and symbolic powers. The specific problem related to this paper is for height $2$ ideals in the standard graded polynomial ring $k[x_0,x_1,x_2]$, namely those ideals that define points in $\PP^2$.
Craig Huneke asked if the containment
\begin{equation}\label{eq: THE containment}
    I^{(3)} \subseteq I^{2}
\end{equation}
always holds, where $I^{(3)}$ denotes the {\it third symbolic power of $I$}. Henceforth \Cref{eq: THE containment} is referred to as {\it the containment problem.} While \Cref{eq: THE containment} holds for general points in $\PP^2$ \cite[Remark 4.3]{bocciharbourne}, the first counterexample to this containment was given in \cite{DST2013}. Sets of points that fail \Cref{eq: THE containment} are special and demonstrate high levels of symmetry. The containment problem, and its variants such as \cite[Conjecture 4.1.1]{HarbourneHuneke} and \cite[Conjecture 2.1]{GRIFOstableharbourne}, have generated a plethora of research; see \cite{AlexandraBen, AlexandraBrianUwe, DUMNICKIResurgence, CounterExampleREALS, AlexandraBrian, justynaMalara, drabkincolinearcodim2, AKESSEH201744}. A survey of the results on the containment problem can be found in \cite{TomasJustyna}.

Both real and projective B\"or\"oczky configurations exhibit extremal behavior in algebraic and combinatorial contexts.
Let $\mathcal{B}_n$ denote the {\it set of triple points} of $B_n$ with $I_n$ as its homogenous defining ideal. The main observation of \cite{CounterExampleREALS} is the fact that $I_{12}^{(3)} \not \subseteq I_{12}^2$; this was the first counterexample to the containment problem over the real numbers. It is noteworthy that Green and Tao's interest in the projective B\"or\"oczky is its extreme behavior in the context of the Dirac-Motzkin conjecture \cite[p.10]{greentao}:
\begin{quote}
    These (B\"or\"oczky) sets, it turns out, provide the example of non-colinear sets of $n$ points with fewest number of ordinary lines, at least for large $n$.
\end{quote}
In fact they are the only configurations, up to projective equivalence, that achieve the lower bound for the number of ordinary lines, which contain ordinary points only, of a configuration of lines and point as conjectured by Dirac and Motzkin; see \cite[Theorem 2.2]{greentao}. 

We summarize the known results about the B\"or\"oczky configuration. The authors in \cite{farnik2017containmentproblemcombinatorics} consider $\mathcal{B}_{12}$ and another configuration of lines and points to show combinatorial data alone does not determine the containment $I^{(3)}\subseteq I^2$. For small values of $n$, the authors in \cite{BoroczkyParameterSpaceOG, BoroczkyParameterSpace2} consider a parameter space of $\mathcal{B}_n$, in order to investigate other rational counterexamples to the containment problem, and \cite{JakubFreeness} studies some combinatorial aspects of B\"or\"oczky configuration. Finally, \cite{Kabatupto12} shows, albeit using Macaulay2, that $I_n$ satisfies \Cref{eq: THE containment} for $1\leq n \leq 11$. 

This paper aims to start the first systematic study of the set $\mathcal{B}_n$ and the algebraic properties of $I_n$ for {\it arbitrary values of} $n$. The following is a summary of our results. Assuming that $3 \mid n$, \Cref{prop: D6 action on triple pts} shows the set of triple intersection points $\mathcal{B}_n$ admits an action of $D_6$ and \Cref{prop: num orbits Bn} computes the number of orbits of size $3$ and $6$ of this action. Next, \Cref{Waldschmit Const n/3} computes the Waldschmit constant of $I_n$. The Waldschmit constant is an important algebraic invariant of an ideal which yields a natural lower bound for the degree of the minimal generators. In \Cref{cor: bounds on In} we use the properties of the weighted projective space $\PP(1,2,3)$ to give an upper bound for the degree of the minimal generators of $I_n$. \Cref{thm: prod lines Unique} shows the smallest degree form in $I_n^{(3)}$ is the product of the lines of $B_n$, the vector space $(I_n^{(3)})_n$ is one dimensional, generated by this product and admits the alternating sign representation of $D_6$. In \Cref{sec: elliptified Boroczky} we describe how applying this construction to points on an elliptic curve resulted in a new counterexample to \Cref{eq: THE containment}. Finally, \Cref{sec: open problems} presents some open problems.

\begin{acknowledgments}
The authors thank Alexandra Seceleanu for comments on earlier versions of this paper. We thank Tomasz Szemberg for helpful discussions during the UWE Fest conference. Roshan Zamir received support from the NSF grants DMS–2401482 and DMS-2342256 RTG: Commutative Algebra at Nebraska.     
\end{acknowledgments}

\section{Action of $D_6$ on $\mathbb{P}^2$} \label{sec: Action}
Fix $k$ to be the field of complex numbers. Consider the natural action of $D_6$, the dihedral group of order $6$, on $\mathbb{A}_k^2$ given by the symmetries of an equilateral triangle. Fix an embedding of $\mathbb{A}^2$ into $\mathbb{P}^2$ via 
\begin{equation}\label{eq: A2 to P2}
     (a,b) \to [a:b:1]. 
\end{equation}
The induced action of $D_6$ on $\mathbb{P}^2$ is given by $g\cdot[a:b:c]=[g(a,b):c]$ and is referred to as  ``the action of $D_6$ on $\PP^2$.'' Equivalently, one can identify $D_6$ with a subgroup of PGL$_3(k)$ generated by 
\begin{equation*}
    \left \langle \left [ \begin{matrix}
        \cos(\frac{2\pi}{3}) & -\sin(\frac{2\pi}{3})& 0 \\
        \sin(\frac{2\pi}{3}) & \cos(\frac{2\pi}{3})& 0 \\
        0 & 0& 1 
    \end{matrix} \right ],  \left [ \begin{matrix}
        1 & 0& 0 \\
        0 & -1& 0 \\
        0 & 0& 1 
    \end{matrix} \right ] \right \rangle
\end{equation*}
This induces an action of $D_6$ on $U_2$, the standard affine chart of $\mathbb{P}^2$ identified with $\AA^2$. 
The induced action on $R=k[x,y,z]$, the homogeneous coordinate ring of $\mathbb{P}^2$, is given by 
\begin{equation}\label{eq: action on poly}
    \begin{split}
        D_6 \times R &\to R\\
        \left (g , \left [ \begin{matrix}
            x\\
            y\\
            z
        \end{matrix}\right ] \right ) &\to g^{-1}\cdot 
        \left [\begin{matrix}
            x\\
            y\\
            z
        \end{matrix} \right ]
    \end{split}
\end{equation}
Note the equations $y=0$, $y-\sqrt{3}x=0$ and $y+\sqrt{3}x=0$ define the three reflection lines of $D_6$, denoted by $l_0$, $l_1$ and $l_2$ respectively, and the $l_i$'s form an orbit under the action defined in \Cref{eq: action on poly}. 
The following notation is in place unless otherwise stated:
\begin{equation}\label{eq: Gens Inv Ring}
\begin{split}    
    p&=\prod_{i=1}^3 l_i=y(y-\sqrt{3}x)(y+\sqrt{3}x)\\
    u&=x^2+y^2\\
    v&=x(x-\sqrt{3}y)(x+\sqrt{3}y)
\end{split}    
\end{equation}
\begin{remark}\label{rmk: p generates skw-inv}
    An element $f \in R$ is a {\it skew-invariant element of the action of $D_6$} if for all $g \in D_6$ either $ g\cdot f = f \text{ or } g\cdot f=-f$, see \cite[p.214]{kane01Reflection}. One can show that $p$ is skew-invariant and if $h$ in $R$ is also skew-invariant, then $p$  divides $h$, see \cite[Proposition B, p.216]{kane01Reflection}.
\end{remark}
Observe $u$ and $v$ in \Cref{eq: Gens Inv Ring} are the defining equation of a circle and the product of lines perpendicular to $l_i$'s respectively. Let $S$ denote $R^{D_6}$, the ring of invariants of the action.
\bpr \label{prop: ring of invarints}
The ring of invariants $S$ is equal to $k[z,u,v]$ and $S$ is regular ring.
\epr 
\bpf 
The proof of the equality $S=k[z,u,v]$ can be found in \cite[Example 2, p.93]{smithinvpoly}. By the Shephard-Todd-Chevalley theorem \cite[Theorem 7.4.1]{smithinvpoly} $S$ is a polynomial ring.
\epf 

\bpr \label{prop: Orbitsize}
An orbit $O$ of this action of $D_6$ has cardinality $1,3$ or $6$. 
\epr 

\bpf
By the Orbit-Stabilizer theorem $|O|$ divides $6$.
If $p\neq [0:0:1]$ and $p\in O$, the rotations of $D_6$ send $p$ to two distinct points thus $3\leq |O|$.
\epf

\begin{remark}\label{rmk: Action Facts}
The following observations are straightforward to check.
\begin{enumerate}
\item Since the action fixes the hyperplane defined by $z$, identify $O$ with its image in $U_2 \cong \AA^2$. If $|O|=3$ or $6$ then elements of $O$ lie on a circle because elements of $D_6$ are isometries which preserve distance. Therefore all elements of an orbit have the same distance from origin.
\item Suppose $O$ is an orbit of size three. Then $O$ has a stabilizer of size $2$, therefore it is fixed by some reflection. This implies that $O$ lies on $V(p)$, the vanishing set of $p$, and any orbit of size $3$ has a point on the $x$-axis. Moreover, each reflection line contains exactly one point of an orbit of size $3$. An orbit of size $6$ can not have a point on $V(p)$ since its stabilizer has size $1$. In particular, $p$ does not vanish on the points of any orbit of size $6$.

\item An invariant polynomial evaluates to the same value on the points of a single orbit. Let $\omega \in k[x,y,z]^{D_6}$, $O$ be an orbit of any size, $a,b \in O$ with $a \neq b$, and $g \in D_6$ such that $g\cdot a=b$. Then
\begin{equation*}
    \omega(a)=\omega(g\cdot b)=g^{-1}\cdot \omega(b)=\omega(b)
\end{equation*}  

\noindent where the last equality follows because $\omega$ is invariant. Therefore any invariant polynomial vanishing on a point of an orbit vanishes on the entire orbit.

\end{enumerate}
\end{remark}
One of our techniques, described in \Cref{obs: P2 to P123}, is realizing the orbits of this action as points in a weighted projective space which we now define.

\begin{definition}[Weighted Projective Space]
    Fix the set $\{a_0,\ldots, a_n \}$ where $a_i \in \Z_{>0}$ for all $i$ and let $Q=k[x_0,\ldots,x_n]$ where $\deg(x_i)=a_i$. The {\it weighted projective space $\PP(a_0,\ldots,a_n)$} is defined as $\Proj(Q)$. 
For a thorough treatment of weighted projective spaces see \cite{hosgoodWPV, DolgachevWPS}.
\end{definition}

\begin{observation}\label{obs: P2 to P123}
Observe $\Proj(S)=\PP(1,2,3)$ is the weighted projective plane with weights $1$, $2$ and $3$.
The map 
\begin{equation*}
    k[z,u,v] \hookrightarrow R
\end{equation*}
induces the map
\begin{equation}\label{eq: P2 to P123}
\begin{split}    
    \phi:\PP^2 &\to \PP(1,2,3),\\
    p=[p_0:p_1:p_2] &\to [z(p):u(p):v(p)]. 
\end{split}
\end{equation}
The map in \Cref{eq: P2 to P123} takes orbits of the action to points of $\PP(1,2,3)$. It is given by invariant polynomials which, by \Cref{rmk: Action Facts}-(3), evaluate to the same value on the points of an orbit.
\end{observation}

The next lemma is a criterion for equality of ideals in a polynomial ring with finitely many variables.
\begin{lemma}[\cite{ENGHETA}]\label{lem: Engheta-Mult Equality}
    Let $R$ be a polynomial ring with a finite number of variables over any field and $J \subset R$ an unmixed ideal. If $I \subset R$ is an ideal containing $J$, with the same multiplicity and height as $J$, then $J=I$.
\end{lemma}

\bpr\label{Single6} 
The defining ideal of an orbit of size $6$ is a complete intersection with two invariant generators of degrees $2$ and $3$.

\epr 

\bpf 
Let $O$ be an orbit of size $6$ with $I$ as its defining ideal. By \Cref{rmk: Action Facts} parts $1$ and $2$, $O$ lies on a circle, which is invariant under the action, defined by $u-r^2z^2$ for some $r >0$. Let $\gamma=v(a)$ for some $a\in O$. Since $v$ is an invariant element, by \cref{rmk: Action Facts}-(3), $\gamma=v(b)$ for all $b\in O$. Furthermore $z(a)=1$ as a result of the map in \Cref{eq: A2 to P2}.
Let $J=(u-r^2z^2, v - \gamma z^3)$. One can check that $J$ is generated by a regular sequence thus $e(R/J)=6$. Clearly $J \subseteq I$ and both ideals have multiplicity $6$. Since height($I$)=height($J$)=2, both ideals define points in $\mathbb{P}^2$ and therefore they are unmixed. The result follows by \Cref{lem: Engheta-Mult Equality}. \epf

One can show the ideal defining an orbit of size $3$ is generated by three quadrics. For example one can choose the pairwise product of the three lines through each pair of the points.
\section{The B\"or\"oczky construction}\label{sec: borocky construction}
\begin{figure}[h!]
    \centering
    \includegraphics[width=0.45\linewidth]{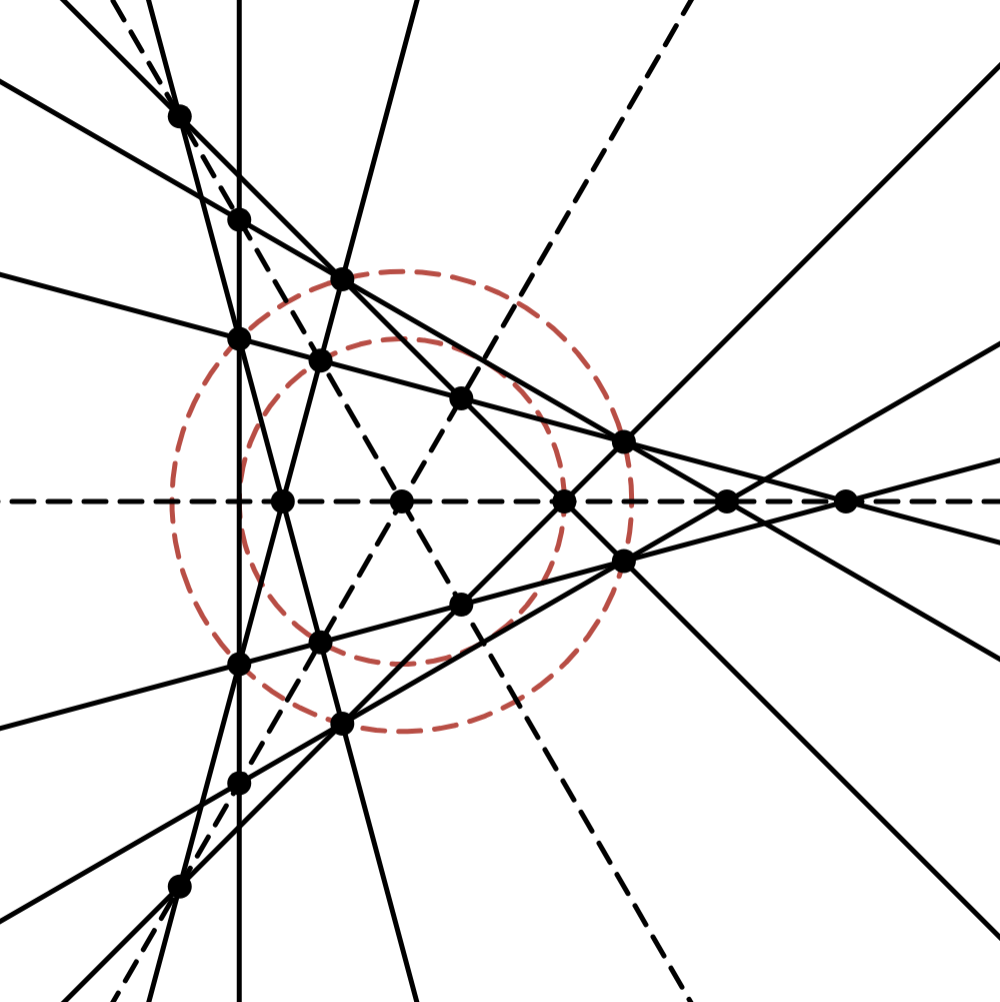}
    \includegraphics[width=0.45\linewidth]{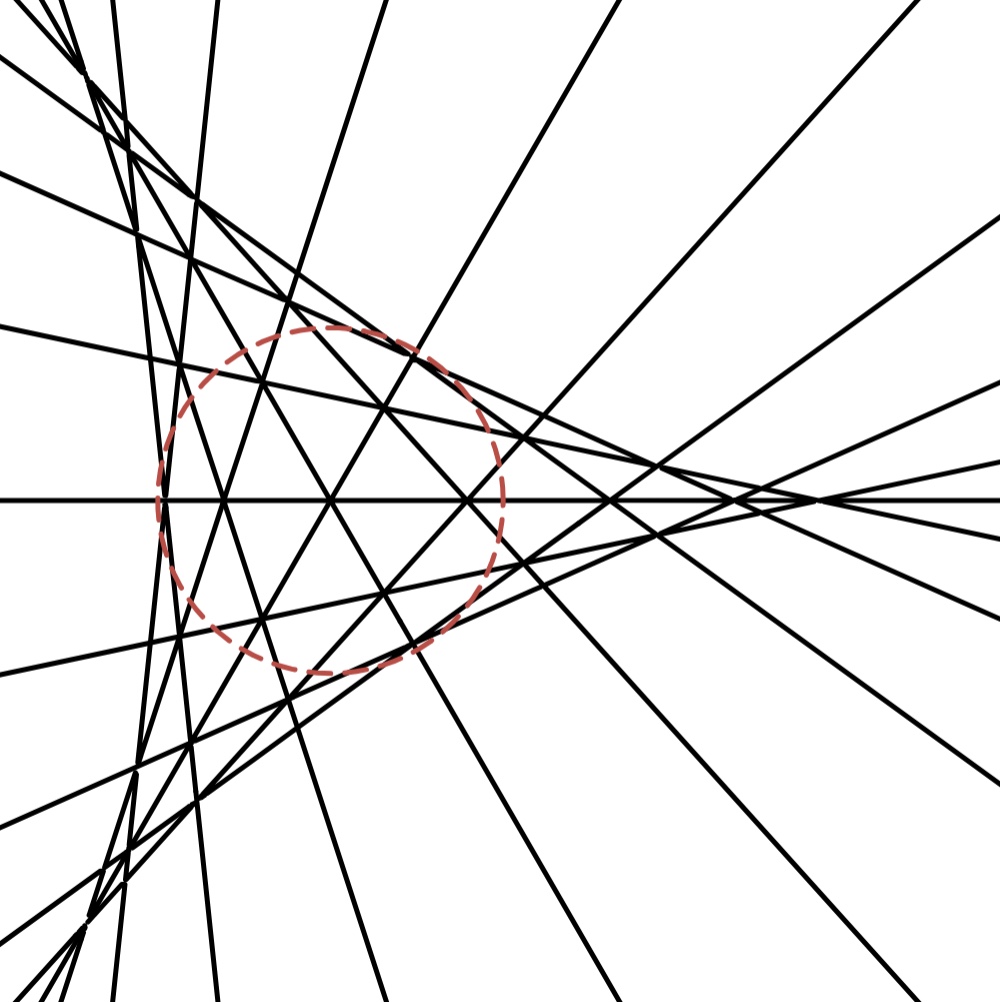}
    \caption{The configurations $B_{12}$(left) and $B_{15}$(right)}
    \label{Fig: B12}
\end{figure}

We recall the construction exactly as described in \cite{PalastiFerudi} over the real plane $\mathbb{R}^2$. Fix a circle $\mathcal{C}$ with center $C$ and a point on the circle $P_O$. For any real number $\alpha$, let $P_\alpha$ denote the point obtained by counterclockwise rotating $P_O$ around $C$, with an angle of $\alpha$. Let $L_\alpha$ denote the real line connecting $P_\alpha$ and $P_{\pi - 2\alpha}$. If $\alpha \equiv \pi-2\alpha \pmod {2\pi}$ then $L_{\alpha}$ is the tangent line to $\mathcal{C}$. Fix $\mathcal{C}$ as the unit circle and $P_O$ as the point $(1,0)$. For any $n\geq 3$ and $0\leq j \leq n-1$, let $L_j=L_{\frac{2\pi j}{n}}$ which connects the points $P_{\frac{2\pi j}{n}}$ and $P_{\pi - \frac{4\pi j}{n}}$. Define the set $\mathcal{L}_n$ as
\begin{equation*}
    \mathcal{L}_n=\{L_{j} \mid 0\leq j\leq n-1 \}.
\end{equation*}
Identify $\R^2$ with the set of complex numbers $\C$ via $(a,b)\to a+bi$ and note for each $i$ the point $P_{\frac{2\pi i}{n}}$ corresponds to $\epsilon=\exp (2\pi i/n)$, an $n$-th root of unity. For $n\geq 3$ where $n$ is even, this construction is a recipe for connecting the $n$ roots of unity because the points $P_{\frac{2\pi j}{n}}$ constitute the the $n$-roots of unity for $0\leq j \leq n-1$. Then each line $L_i$ of $\mathcal{L}_n$ denotes the line $L_{i,\frac{n}{2}-i} $, which connects $P_i$ and $P_{\frac{n}{2}-i}$. For ease of notation we interchangeably use $L_i$ to denote the line and its defining equation. We fix the following notation:
    \begin{itemize}
        \item $B_n$ := Union of lines and all intersection points in the B\"or\"oczky configuration of size $n$.
        \item $\cL_n$ := The set of lines of $B_n$.
        \item $\cB_n$ := The set of triple points of $B_n$.
    \end{itemize}       
\bpr\label{prop: lines meet}
The lines $L_i,L_j,L_k$ of $\mathcal{L}_n$ are concurrent if and only if $n$ divides $i+j+k$. The number of triple points of $B_n$ is given by
\begin{equation}\label{eq: number triple}
    \lfloor n(n-3)/6 \rfloor +1.
\end{equation}
\epr 
\bpf
The first assertion is proved in \cite[Section 2, Lemma]{PalastiFerudi} and \Cref{eq: number triple} follows from \cite[Section 6, Property 4]{PalastiFerudi}.
\epf 
The configurations $B_{12}$ and $B_{15}$, created with Desmos, are demonstrated in \Cref{Fig: B12}. \href{https://www.desmos.com/calculator/jdi7xwvif2}{This link contains our Desmos code}. In $B_{12}$ the dotted lines represent $p$ as defined in \Cref{eq: Gens Inv Ring}; they are part of the $12$ configuration lines. The marked points are the $19$ points of $\mathcal{B}_{12}$.
 
In \Cref{prop: D6 action on triple pts} we show the group action of $D_6$ on $\PP^2$ induces an action on the set $\mathcal{L}_n$ and subsequently $\mathcal{B}_n$. 
\bpr \label{prop: D6 action on triple pts}
Assume $n\geq 3$ and $n$ is divisible by $3$. There is an induced action of $D_6$ on $\mathcal{B}_n$.
\epr 
\bpf 
let $n=3k$ for some $k\geq 1$. Let $\epsilon=\exp(\frac{2\pi i}{2n})$ and $P_r=\epsilon^r$. The set $P'=\{P_{2i} \}_{i=0}^{n-1}$ is in the set of $n$-th roots of unity indexed as part of the $2n$-th roots of unity. The construction applied to the $2n$ roots of unity connects the points $P_{2i}$ and $P_{n-4i}$, where the indices are understood modulo $2n=6k$. Next, we wish to consider the lines of the construction connecting the points in the set $P'$, that is the set of lines
\begin{equation*}
    L'=\{L_{2i,n-4i} \mid 0\leq i \leq n-1 \}.
\end{equation*}
Since $D_6$ is generated by $T_1$, the rotation by $\pi/3$ and $T_2$, the reflection over the $x$-axis, it suffices to show $T_1\cdot L'=L'$ and $T_2\cdot L'=L'$. For any $0\leq i \leq n-1$ we have:
\begin{equation*}
\begin{split}    
    & T_1 \cdot (P_{2i})=P_{2i+2k} \text{ and } T_1\cdot (P_{n-4i})=P_{n-4i+2k}=P_{n-4i-4k}\\
    \Rightarrow &T_1 \cdot (L_{2i,n-4i})=L_{2i+2k,n-4i+2k}=L_{2i+2k,n-4i-4k}=L_{2i+2k,n-2(2i+2k)},\\
    &T_2 \cdot (P_{2i}) = P_{-2i} \text{ and } T_2 \cdot P_{n-4i}=P_{-n+4i}=P_{n-2(-2i)}\\
    \Rightarrow &T_2(L_{2i,n-4i})=L_{-2i,n-2(-2i)}.
\end{split}
\end{equation*}
Let $L'_i=L_{2i,n-4i}$ for any $0\leq i \leq n-1$.
Finally, we observe the points of $\mathcal{B}_n$ are mapped back to themselves since if $P$ is on $L'_r \cap L'_j \cap L'_k$ then $T_i(P)$ $T_i(L'_r) \cap T_i(L'_j) \cap T_i(L'_k)$ for $i=1,2$.
\epf 
A subtlety in the proof of \Cref{prop: D6 action on triple pts} is its equivocation of the B\"or\"oczky construction applied to $n$ roots of unity in and of themselves and as part of the $2n$ roots of unity. These constructions are equivalent although the indexing of the lines differs.

The action of $D_6$ can be observed in \Cref{Fig: B12} where $\mathcal{B}_{12}$ consists of one orbit of size $6$, four orbits of size $3$ and one orbit of size $1$. The dotted circle with larger radius corresponds to the degree $2$ invariant generator from \Cref{Single6}. It seems like a natural question to compute the number of orbits of a given size of $\mathcal{B}_n$. In particular, the number of orbits of size $6$ is fundamental to the proof of \Cref{cor: bounds on In}. Let $\mathcal{O}_6(n)$ and $\mathcal{O}_3(n)$ denote the number of orbits of size $6$ and $3$ of $\mathcal{B}_n$ respectively.

\bpr\label{prop: num orbits Bn}
For any positive integer $n\geq 3$ where $3\mid n$, the functions $\cO_6(n)$ and $\cO_3(n)$ can be computed via 
\begin{align*}
    \cO_3(n)&=\left\lfloor\frac{n-1}{2}\right\rfloor-1,\\
    \cO_6(n)&=\frac{n(n-3)}{36}-\frac{1}{2}\left\lfloor\frac{n-1}{2}\right\rfloor+\frac{1}{2}.\\
\end{align*}
\epr
\bpf
The assumption $3 \mid n$ implies the quantity $\frac{n(n-3)}{6}$ is an integer and \Cref{eq: number triple} implies the set $\cB_n$ consists of $\frac{n(n-3)}{6}+1$ points. By \Cref{prop: D6 action on triple pts} the set of triple points is acted upon by $D_6$ and the point $(0,0)$ forms the only orbit of size $1$.

The orbits of size 3 comprise the points on $v(p)$ where recall $p$ is the product of the lines $y=0$, $\sqrt{3}x+y=0$, and $\sqrt{3}x-y=0$. Without loss of generality consider the triple points along the line $y=0$. Fix an $n$ that satisfies the hypothesis of the lemma. By \Cref{prop: lines meet} every line $L_i:=L_{i,\frac{n}{2}-2i}$ meets $L_{n-i}$ at $y=0$. When $\frac{n}{2}$ is an integer the line $L_{\frac{n}{2}}$ does not have a triple point on the line $y=0$. Thus, out of the $n-1$ total lines of $\cL_n\setminus\{V(y)\}$, $\lfloor\frac{n-1}{2}\rfloor$ have a triple point on $V(y)$. One of these points is the origin $(0,0)$, and so there are a total of $\lfloor\frac{n-1}{2}\rfloor-1$ points of orbits of size 3 on $V(y)$, yielding the total number of orbits of size 3 as 
\begin{equation*}
    \cO_3(n)=\left\lfloor\frac{n-1}{2}\right\rfloor-1.
\end{equation*}
Thus there are $3\cdot\cO_3(n)=3\lfloor\frac{n-1}{2}\rfloor-3$ points in orbits of size $3$. Straightforward calculations yield the number of orbits of size $6$ as:
\begin{equation*}
    \cO_6(n)=\frac{n(n-3)}{36}-\frac{1}{2}\left\lfloor\frac{n-1}{2}\right\rfloor+\frac{1}{2}. \qedhere
\end{equation*}
\epf
\bcor\label{prop: formula orbits Bn} For any positive even integer $n\geq 6$ where $6\mid n$, the functions $\cO_3(n)$ and $\cO_6(n)$ are given by
\begin{equation}\label{prop: formula orbits Bn orbits Bn}
\begin{split}
    \cO_3(n)&=\frac{n}{2}-2,\\
    \cO_6(n)&=\frac{(n-6)^2}{36},
    \end{split}
\end{equation}
\ecor 
We finish this section with two observations about the set $\mathcal{L}_n$. If $6 \mid n$ then $\mathcal{L}_n$ contains exactly three tangent lines to the unit circle. This follows from solving the modular equation
\begin{equation*}
    i \equiv \frac{n}{2}-2i{\pmod n} \text{ for } 0\leq i \leq n-1,
\end{equation*}
where the solutions are $\frac{n}{6},\frac{n}{2},\frac{5n}{6}$. If $3\mid n$ then the lines of $p$ are always contained in $\mathcal{L}_n$; they correspond to lines $L_0$,$L_{\frac{n}{3}},L_{\frac{2n}{3}}$.
\section{Main Results}\label{sec: Borocky General}
Let $I_n\subseteq k[x,y,z]$ be the defining ideal of $\mathcal{B}_n$, the set of triple intersection points of $B_n$. The main results of this section are \Cref{Waldschmit Const n/3}, which uses \Cref{lemma: least number of points} to compute the Waldschmit constant of $I_n$, \Cref{cor: bounds on In} which uses the weighted projective space $\PP(1,2,3)$ to bound the minimal degree of a generator of $I_n$, and \Cref{thm: prod lines Unique} which shows the product of the $n$ configuration lines is the unique element of degree $n$ in $I_n^{(3)}$.

Fix a non-negative integer $N$ and let $k[\PP^N]$ denote the homogeneous coordinate ring of $\PP^N$. Recall for a non-zero, homogeneous ideal $J \subsetneq k[\PP^N]$, the {\it Waldschmit constant of $J$}, denoted by $\widehat{\alpha}(J)$, is defined as the limit
\begin{equation*}
    \widehat{\alpha}(J)=\lim_{m\to \infty}\cfrac{\alpha(J^{(m)})}{m}=\inf_{m\geq 0}\cfrac{\alpha(J^{(m)})}{m},
\end{equation*}
where $\alpha(J^{(m)})$ represents the smallest degree of a non-zero element of $J^{(m)}$, the $m$-th symbolic power of $J$. The existence of the limit and the second equality follow from \cite[Lemma 2.3.1]{bocciharbourne}. 
\blm \label{lemma: least number of points}
Let $n\geq 8$ be an even integer. Each line of $\cL_n$ contains at least $\LF \frac{n}{3}\RF $ configuration points.
\elm
\bpf 
By \Cref{prop: lines meet} three lines $L_i$, $L_j$ and $L_k$ are concurrent if and only if $i+j+k \equiv 0 \pmod n$. Fix $2\leq i$ and let $|L_i|$ denote the number of points of $\mathcal{B}_n$ on $L_i$. Since $ 0\leq i,j,k \leq n-1$, computing the quantity $|L_i|$ is equivalent to counting the number of distinct non-negative integer solutions to the equations
\begin{equation} \label{eq: n-i}
    j+k=n-i 
\end{equation}
\begin{equation} \label{eq: 2n-i}
    j+k=2n-i
\end{equation}
where $i,j,k$ are pairwise distinct as no line $L_i$ intersects itself at a configuration point. Moreover, $(j,k)$ and $(k,j)$ represent the same intersection. The assumption that $n$ is even implies $n-i$ and $2n-i$ are even if and only if $i$ is even. We first compute the number of solutions to \Cref{eq: n-i}:
\begin{equation}\label{eq: solutions n-i}
\begin{split}    
    &i \text{ is even: }j \in \left [0,\frac{n-i}{2} \right ] \text{ and } k \in \left [\frac{n-i}{2},n-i \right ] \Rightarrow \frac{n-i}{2} \text{ solutions} \\
    &i \text{ is odd: } j \in \left [0,\frac{n-i-1}{2} \right ] \text{ and } k \in \left [\frac{n-i+1}{2},n-i\right ] \Rightarrow  \frac{n-i+1}{2}  \text{ solutions} 
\end{split}
\end{equation}
Since $j$ and $k$ must be distinct, the solution $\left ( \frac{n-i}{2}, \frac{n-i}{2} \right )$ is excluded from the case $i$ is even.

Similarly, to count the number of solutions to \Cref{eq: 2n-i} one only needs to consider the range $0\leq j \leq \LF \frac{2n-i}{2} \RF$.
Next, we claim that if $(j,k)$ is a solution to \Cref{eq: 2n-i}, then $ n-i+1 \leq j$. Assume to the contrary that $ j < n-i+1$ and observe the following contradicts that $k \leq n-1$.
\begin{equation*}
    j < n-i+1 \iff j < 2n-i-n+1 \iff   n-1< 2n-i-j\underset{\ref{eq: 2n-i}}{=}k.
\end{equation*}
The claim follows. In \Cref{eq: solutions 2n-i} we compute the number of solutions to \Cref{eq: 2n-i}. When $i$ is even, the claim implies $j \in \left [n-i+1, n-\frac{i}{2} \right ]$ and the the pair $(n-\frac{i}{2},n-\frac{i}{2})$ is excluded for distinctness.
\begin{equation}\label{eq: solutions 2n-i}
\begin{split}
    &i \text{ is even: }j \in \left [n-i+1, n-\frac{i}{2} \right ] \text{ and } k \in \left [n-\frac{i}{2}, n-1 \right ]  \Rightarrow \frac{i}{2}-1 \text{ solutions} \\
    &i \text{ is odd: }j \in \left [n-i+1, {\textstyle\frac{2n-i-1}{2}}\right ] \text{ and } k \in \left [{\textstyle\frac{2n-i+1}{2}},n-1 \right ] \Rightarrow \frac{i-1}{2}\text{ solutions}    
\end{split}
\end{equation}
Thus far there are a total of $\frac{n}{2}-1$ solutions when $i$ is even and $\frac{n}{2}$ solutions when $i$ is odd. As $i,j,k$ are pairwise distinct, we need to count and subtract the cases when $i=j$ or $i=k$. For a fixed $i$, the intervals in \Cref{eq: solutions n-i} and \Cref{eq: solutions 2n-i} are disjoint and they cover all integers from $0$ to $n-1$, thus $i$ will equal $j$ or $k$ for some solution pair $(j,k)$ and we must subtract one more from each solution count. However, if:
\begin{itemize}
    \item $\frac{n-i}{2}=i$ which is true if and only if $i$ is even and $i=\frac{n}{3}$, or
    \item $n-\frac{i}{2}=i$ which is true if and only if $i$ is even and $i=\frac{2n}{3}$,
\end{itemize}
then we have already accounted for the repeated solution and we don't need to subtract anything. For $i=0,1$, \Cref{eq: 2n-i} has no solutions and a simple counting shows \Cref{eq: n-i} has $\frac{n}{2}-1$ solutions.  
To summarize, if $i\geq 2$ and $i$ is even, we have
\begin{equation*}
     |L_i| =
    \begin{cases}
    \frac{n}{2}-1 & \text{ if } i=\frac{n}{3} \text{ or } i=\frac{2}{3}n, \\
    \frac{n}{2}-2 &\text{ otherwise .}
\end{cases}
\end{equation*}
And if $i\in \{0,1\}$ or $i$ is an odd number we have
\begin{equation*}
    |L_i| = \frac{n}{2}-1.
\end{equation*}
Finally, the proof concludes since $\frac{n}{2}-2\leq |L_i|$ for all cases of $i$ and 
\begin{equation}\label{eq: pts on Ln}
    \LF \frac{n}{3}\RF \leq \frac{n}{2}-2 \text{ when } n\geq 8. \qedhere
\end{equation}
\epf 

\bcor \label{cor: points of Ln}
For $n\geq 14$ where $n$ is even, each line of $\mathcal{L}_n$ contains at least $\LF \frac{n}{3} \RF+1$ configuration points.
\ecor

\bpf 
The hypothesis implies \Cref{eq: pts on Ln} is a strict inequality.
\epf 

\bthm \label{Waldschmit Const n/3}
For any $n\geq 12$ where $6 \mid n$, we have $
    \widehat{\alpha}(I_n)= \frac{n}{3}.$
\ethm 
\bpf 
First we establish the inequality $\widehat{\alpha}(I_n)\leq \frac{n}{3}$. Observe the product $\prod_{i=1}^n L_i$ is in $I_n^{(3)}$ thus $\alpha(I_n^{(3)})\leq n$ and \Cref{eq: ineq implic} yields $\widehat{\alpha}(I_n)\leq \frac{n}{3}$.
\begin{equation}\label{eq: ineq implic}
    \alpha(I_n^{(3m)})\leq m\alpha(I_n^{(3)})\leq mn\Rightarrow  \widehat{\alpha}(I_n)=\lim_{m\to \infty} \frac{ \alpha(I_n^{(3m)})}{3m}\leq \frac{n}{3}.
\end{equation}
Assume for the sake of contradiction that $\widehat{\alpha}(I_n) <  \frac{n}{3}$. By definition of the infimum there exists an $m$ such that
\begin{equation*}
\frac{\alpha(I_n^{(m)} )}{m} < \frac{n}{3}\Rightarrow
\alpha(I_n^{(m)} ) \leq m \frac{n}{3}   -1.
\end{equation*}
Let $G \in I_n^{(m)}$ such that $\deg(G)\leq m\frac{n}{3}-1$ and note $G$ vanishes on $\mathcal{B}_n$ to order at least $m$. Recall $L_i$ denotes both the line $L_{i,\frac{n}{2}-i}$ and its defining equation. By \Cref{lemma: least number of points}, $G$ intersects $L_i$ at at least $\frac{n}{3}$ points thus $e(L_i \cap G)\geq  \frac{n}{3} m $, where $e$ denotes the intersection multiplicity. By Bezout's theorem $L_i \mid G$ for all $i$. Let $G'=G\backslash (\prod_{i=1}^n L_i)$ which implies $\deg(G')=\deg(G)-n$. Since $G'$ vanishes on $\mathcal{B}_n$ to order at least $m-3$ we get that $G' \in I_n^{(m-3)}$. Thus 
\begin{equation*}
    \alpha(I_n^{(m-3)} )  \leq \deg(G)-n\leq m \frac{n}{3}  -1-n = \frac{n}{3}(m-3)-1
\end{equation*}
Continuing in this fashion we arrive at the cases
\begin{equation}\label{eq: upper bound F}
    m-3 = \begin{cases}
        1 \\
        2 \\
        3
    \end{cases}
    \Rightarrow \deg(F) \leq 
    \begin{cases}
        \frac{n}{3}-1\\
        \frac{2n}{3}-1\\
        n-1
    \end{cases}
\end{equation}
where $F$ is an element of minimal degree in $I_n^{(m-3)}$. Let $\gamma \in \{1,2,3 \}$. In each instance of $\gamma$, \Cref{lemma: least number of points} implies $e(F \cap L_i) \geq \frac{n}{3}(\gamma)$ and thus $F$ is divisible by all $n$ lines of the configuration. Thus $\deg(F)\geq n$ which contradicts 
\Cref{eq: upper bound F}. Hence $\widehat{\alpha}(I_n) = \frac{n}{3}$. 
\epf
Observe that $\widehat{\alpha}(I_n)$ is a natural lower bound for $\alpha(I_n)$. We introduce enough preliminary results to give an upper bound for $\alpha(I_n)$ in \Cref{cor: bounds on In}.
\begin{remark}\label{rmk: S_d}
    
Let $S=k[z,u,v]$. Recall from \Cref{eq: Gens Inv Ring} that the degrees of $z$, $u$ and $v$ are $1,2$, and $3$ respectively. Define $s_d$ to be the Hilbert function of $S$ in degree $d$, that is $s_d:=\dim_k((S)_d)$. For any $d \in \Z_{\geq 0}$, $s_d$ represents the number of linearly independent monomials of degree $d$ in the non-standard graded ring $S$. By \cite[p.547]{stanley1}, $s_d$ can be computed via the formula:
\begin{equation}\label{eq: s_d}
    s_d = \LF \frac{d^2}{12} + \frac{d}{2} + 1 \RF.
\end{equation}
\end{remark}
\bpr \label{prop: smallest inv}
The ideal $I_n$ contains a skew-invariant element of degree $d+3$ where $d$ is the smallest integer satisfying
\begin{equation}\label{eq: sd > pts}
   s_d-\frac{(n-6)^2}{36} \geq 1 \iff  \LF\frac{d^2}{12} + \frac{d}{2} \RF \geq \frac{(n-6)^2}{36}.
\end{equation}
\epr  

\bpf 
By \Cref{prop: formula orbits Bn orbits Bn}, the set $\cB_n$ contains $\frac{(n-6)^2}{36}$ orbits of size $6$. By \Cref{rmk: Action Facts}-(3), invariant polynomials evaluate to the same value on the points of the same orbit. Therefore, applying the map in \Cref{eq: P2 to P123}, which is given by invariant polynomials, to each orbit of size $6$ results in $r=\frac{(n-6)^2}{36}$ points in $\PP(1,2,3)$, denoted by $X=\{p_i \}_{i=1}^r \subseteq \PP(1,2,3)$. Let $I_X \subseteq S$ be the graded homogeneous ideal of $X$. Clearly $(I_X)_d$, the $d$-th graded piece of $I_X$, is equal to the kernel of the evaluation map
\begin{equation*}
\begin{split}    
    \ev_X (d): &(S)_d \to k^r,\\
                &f \to (f(p_i))_{i=1}^r.
\end{split}
\end{equation*}
The Rank-Nullity theorem implies the inequality 
\begin{equation*}
    \max \{s_d - r, 0\} \leq \dim((I_X)_d).
\end{equation*}
Since $n$ is fixed and \Cref{eq: s_d} shows $s_d$ is increasing, there exists an integer $d$ satisfying \Cref{eq: sd > pts}. By the well-ordering principle one can choose the smallest such integer $d$. Therefore there exists an invariant form, $F$, of degree $d$ that vanishes on the orbits of size $6$ of $\cB_n$. By \Cref{rmk: Action Facts}-(2), $p$ as defined in \Cref{eq: Gens Inv Ring} vanishes on all orbits of size $3$ and the origin. Thus $p\cdot F$ is the desired element. 
\epf 

\bcor\label{cor: bounds on In}
Let $d$ be as defined in \Cref{prop: smallest inv}. Then
\begin{equation}\label{eq: bounds on min gen}
 \frac{n}{3}\leq \alpha(I_n) \leq d+3.
\end{equation}
\ecor

\bpf
The inequality $\frac{n}{3}\leq \alpha(I_n)$ follows because $I_n=I_n^{(1)}$ and the Waldschmit constant is a lower bound for $\alpha(I_n^{(1)})$. The inequality $\alpha(I_n) \leq d+3$ follows from \Cref{prop: smallest inv}.
\epf 
\begin{remark}
    Observe the upper bound in \Cref{eq: bounds on min gen} is sharper than the bound one would get by viewing the orbits of size $6$ as general points in $\PP^2$, finding a curve through them and taking the union with $p$. This is because the function in \Cref{eq: s_d} grows slower than $\binom{d+2}{d}$. For example, when $n=12$, one can use \Cref{eq: s_d} to show \Cref{eq: bounds on min gen} yields 
    \begin{equation}\label{eq: bnds I12}
        4\leq \alpha(I_{12}) \leq 5.
    \end{equation}
    Note $6$ general points in $\PP^2$ do not necessarily define a conic, but they define a pencil of cubics. Fix a cubic $C$. Then $p\cdot C$ would vanish on $\cB_{12}$ but $\deg(p\cdot C)=6$. In fact, one can show $\alpha(I_{12})=5$. Assume for the sake of contradiction that $\alpha(I_{12})=4$ and let $F$ be a generator of minimal degree. Note the lines defining $p$ are the three lines through the origin in \Cref{Fig: B12}, each of which contains $5$ configuration points. Bezout's theorem implies each line of $p$ divides $F$ and we get $F=p\cdot l$ for some linear form $l$. By \Cref{rmk: Action Facts}-(2), $p$ does not vanish on the points of an orbit of size $6$, thus $l$ must vanish on the orbit of size $6$ in $\cB_{12}$. But an orbit of size $6$ cannot be on a line by \Cref{rmk: Action Facts}-(1). Alternatively, \Cref{Single6} shows the ideal defining an orbit of size $6$ is generated in degrees $2$ and $3$ and contains no linear forms.
\end{remark}
Computations with Macaulay2 show for any $n\geq 12$ the product of the $n$ lines of the configuration is the {\it unique} element in $(I_n)^{(3)}$ that is not in $I_n^2$. \Cref{thm: prod lines Unique} shows this product is actually the smallest degree form vanishing to order $3$ on the points of $\mathcal{B}_n$.
\bpr\label{thm: prod lines Unique}
For any $n\geq 12$ where $6 \mid n$, we have $\alpha(I_n^{(3)})=n$ and $(I_n^{(3)})_n$ is a one dimensional vector space generated by the product of all the configuration lines i.e. 
\begin{equation}\label{eq: In symb 3 deg n}
    (I_n^{(3)})_n=\Span_k\left \{ \prod_{i=1}^n L_i \right \}.    
\end{equation}
Furthermore, $(I_n^{(3)})_n$ admits the alternating-sign representation of $D_6$.
\epr

\bpf 
The inequality $\frac{n}{3} \leq \frac{\alpha(I_{n}^{(3)})}{3} $ follows from \Cref{Waldschmit Const n/3}, thus $n\leq \alpha(I_{n}^{(3)})$. By the construction outlined in \Cref{sec: Action} the product of the $n$ lines of the configuration is in $I_n^{(3)}$ therefore $\alpha(I_{n}^{(3)})\leq n$ and the inclusion $\supseteq$ in \Cref{eq: In symb 3 deg n} follows. 
Let $n\geq 18$ and $F \in ((I_n)^{(3)})_n$. By \Cref{cor: points of Ln}, each $L_i$ contains at least $\frac{n}{3}+1$ triple points and since $F$ vanishes to order at least $3$ on each triple point we get
\begin{equation*}
    |F \cap L_i|\geq n+3 > n=\deg(F) \text{ for all } i.
\end{equation*}
Therefore Bezout's theorem $L_i \mid F$ for $1\leq i \leq n$ and since $\deg(F)=n$ we get $F=c\cdot \prod_{i=1}^{n} L_i$ for some constant $c$. 
Suppose $n=12$ and let $G\in (I_{12}^{(3)})_{12}$. Observe in \Cref{Fig: B12} that any configuration line $L$ that is {\it not} tangent to the unit circle contains $5$ points of $\mathcal{B}_{12}$ thus $|G\cap L|\geq 15$ and $L\mid G$ by Bezout's theorem. Therefore $G=L'\cdot C$ where $L'$ is the product of the 9 lines not tangent to the unit circle and $C$ is a cubic form. Finally, observe in \Cref{Fig: B12} that $L'$ covers every point of $\mathcal{B}_{12}$ to order $3$ except those located on the tangent lines. Thus, each tangent line intersects $C$ at 4 points and by Bezout's theorem the result follows. Recall the product of the lines $L_i$ indexed by $i=0,\frac{n}{3},\frac{2n}{3}$ is $p$. Therefore $p\mid \prod_{i=1}^{n} L_i$ and \Cref{rmk: p generates skw-inv} implies $\prod_{i=1}^{n} L_i$ is skew-invariant. Finally, a vector space generated by a skew-invariant element admits the alternating sign representation.
\epf 
\begin{remark}\label{rmk: data min degrees gens}
    The Macaulay2 generated data in \Cref{tab: min gen data} reports the degree and the number of the minimal generators for some examples of $I_n$. As far as we know, these computations are new for $n\geq 18$. The number of repeated entries is the number of generators of a particular degree.
    \begin{table}[h!]
    \centering
    \begin{tabular}{|c|c|}
    \hline
    $n$ & Minimal Generators of $I_n$  \\ \hline   
        $12$  & $5,5,5$\\ \hline
           $18$ & $8,8,8,9$ \\ \hline
        $24$ & $11,11,11,12,13$\\ \hline
        $30$ & $14,14,14,15,16,17$\\ \hline
        $36$ & $11$ generators of degree $19$ \\ \hline
    \end{tabular}          
    \caption{Cases of Minimal Generators for $I_n$}
    \label{tab: min gen data}
\end{table}    
\end{remark}

\section{Elliptified B\"or\"oczky}\label{sec: elliptified Boroczky}
In this section we briefly explain how applying the B\"or\"oczky construction to an elliptic curve resulted in a counterexample to $I^{(3)}\subseteq I^2$, which is new to our knowledge. We call this the {\it Elliptified B\"or\"oczky} configuration, constructed over $\mathbb{P}_{\C}^2$ with homogeneous coordinate ring $R=\C[x_0,x_1,x_2]$. Consider the elliptic curve in $\PP^2$ defined by $x^3+y^3+z^3=0$ equipped with the group law for a fixed flex point. One can show under this operation the subgroup generated by all elements of order $n$, denoted $E[n]$, is isomorphic to $\mathbb{Z}_n\times \Z_n$ \cite[Theorem 4.6]{Hartshorne}. We set the notation for this section.
\begin{enumerate}
    \item Let $t$ be a primitive third root of unity i.e. a solution to $t^2+t+1=0$.
    \item Define $s=-\sqrt[3]{2}$ i.e. a solution to $s^3+2=0$.
    \item Let $\alpha=[1:s:1]$ and $\beta=[s:1:t]$.
\end{enumerate} 
When $n=6$, one can show the field extension $\Q(s,t)$ contains the coordinates of all $36$ points of $E[6]$ and $E[6]=\langle \alpha ,\beta \rangle $ with the operation table: 
\small
\begin{table}[h!]
\centering
    \begin{tabular}{|c||c|c|c|c|c|c|}
    \hline
        + &0&$\alpha$&$2\alpha$& $3\alpha$& $4\alpha$& $5\alpha$ \\
        \hline \hline
        0 &$[1:-1:0]$&$[1:s:1]$&$[0:1:-1]$& $[1:1:s]$& $[1:0:-1]$& $[s:1:1]$ \\
        \hline
         $\beta$ &$[s:1:t]$&$[t^2:1:s]$&$[t:s:1]$& $[s:t^2:1]$& $[1:t:s]$& $[1:s:t^2]$ \\
        \hline
         $2\beta$ &$[-1:0:t]$&$[s:t^2:t]$&$[-t:1:0]$& $[t:s:t^2]$& $[0:1:-t^2]$& $[t^2:t:s]$ \\
        \hline
         $3\beta$ &$[t^2:t^2:s]$&$[t:s:t]$&$[s:t^2:t^2]$& $[t:t:s]$& $[t^2:s:t^2]$& $[s:t:t]$ \\
        \hline
         $4\beta$ &$[0:1:-t]$&$[t:t^2:s]$&$[-1:0:t^2]$& $[s:t:t^2]$& $[-t^2:1:0]$& $[t^2:s:t]$ \\
        \hline
         $5\beta$ &$[1:s:t]$&$[s:1:t^2]$&$[t:1:s]$& $[t^2:s:1]$& $[s:t:1]$& $[1:t^2:s]$ \\
        \hline
    \end{tabular}
    \caption{Operation table for $E[6]$}
    \label{tab: Mult Table}
\end{table}  
\normalsize

For $1\leq i_1,i_2,j_1,j_2 \leq 6$, define the line $\ell(i_1,j_1,i_2,j_2)$ as
\begin{equation*}
    \begin{cases}
    \text{the line through $i_1\alpha+j_1\beta$ and $i_2\alpha+j_2\beta$}&i_1\alpha+j_1\beta\neq i_2\alpha+j_2\beta\\
    \text{the tangent line to $E$ at $i_1\alpha+j_1\beta$}&i_1\alpha+j_1\beta= i_2\alpha+j_2\beta\\
\end{cases}
\end{equation*} and 
\begin{equation*}
    L(i,j):=\ell(i,j,3-2i,3-2j)
\end{equation*}
where the entries of $\ell$ are understood modulo $6$.
Consider the configuration of 36 lines 
\begin{equation*}
    \mathcal{L}=\{L(i,j):0\leq i,j\leq 5\}.
\end{equation*}
The lines $L(1,1)$, $L(1,3)$, $L(1,5)$, $L(3,1)$, $L(3,3)$, $L(3,5)$, $L(5,1)$, $L(5,3)$, and $L(5,5)$ are the respective tangent lines to the elliptic curve through the points $[t^2:1:s]$, $[s:t^2:1]$, $[1:s:t^2]$, $[t:s:t]$, $[t:t:s]$, $[s:t:t]$, $[s:1:t^2]$, $[t^2:s:1]$, and $[1:t^2:s]$. Macaulay2 computations in \Cref{section: M2 code} show:
\begin{enumerate}
    \item Up to constant multiples the set $\mathcal{L}$ contains $18$ unique lines.
    \item The lines of $\mathcal{L}$ intersect at $57$ points (unique up to constant multiples) with multiplicity at least $2$, $48$ points with multiplicity exactly $3$. Unlike many other counterexamples to $I^{(3)} \subseteq I^2$ this configuration contains no quadruple, or higher order, intersection points. Let $\mathcal{B}_{E6}$ denoted the set of triple intersection points.
    \item No point of $\mathcal{B}_{E6}$ is located on the elliptic curve $V(x^3+y^3+z^3)$.
    \item Consider the action $S_3 \curvearrowright \mathbb{P}^2$ given by permuting the coordinates i.e. $\sigma \cdot [p_0:p_1:p_2]=[p_{\sigma(0)}: p_{\sigma(1)}:p_{\sigma(2)}]$. One can observe from \Cref{tab: Mult Table} that the restriction of this action maps $E[6]$ back to itself. And one can show both of the sets $\mathcal{L}$ and $\mathcal{B}_{E6}$ are closed under this action of $S_3$. The $48$ points of $E[6]$ are given by the closure of the $S_3$-action on the following $13$ points which form one orbit of size one, one of size two, seven of size three, and four of size six respectively.    
    \begin{table}[h!]
    \centering
    \begin{tabular}{|c|c|}
    \hline
    Size of the orbit & Representative for distinct orbits \\ \hline   
        $1$  & $[1:1:1]$\\ \hline
        $2$ & $[1:t:t^2]$ \\ \hline
        $3$ & \vtop{\hbox{\strut $[1:0:0],[1:1:t],[1:t:t],[1:1:-s^2-t],$}\hbox{\strut $[1:1:-s^2t-1],[1:1:-s^2t^2-t^2],[1:s^2t:s^2t]$}}\\ \hline
        $6$ & \vtop{\hbox{\strut $[1:t:-s^2+t+1],[1:t:-s^2t-t]$,}\hbox{\strut $[1:t:s^2t+s^2-1],[2:2t:-st]$}}\\ \hline
    \end{tabular}          
    \caption{Representatives for distinct orbits of $E[6]$}
    \label{tab:my_label}
\end{table}    
\end{enumerate}
Let $I_{E6}$ be the homogeneous defining ideal of $\mathcal{B}_{E6}$ inside the polynomial ring $\mathbb{Q}(s,t)[x_0,x_1,x_2]$. The Macaulay2 code in \Cref{section: M2 code} shows that
\begin{equation*}
    I_{E6}^{(3)} \not \subseteq I_{E6}^2.
\end{equation*}
In particular, the product of the $18$ unique lines of $\mathcal{L}$ is again the element contained in $I_{E6}^{(3)}$ and not in $I_{E6}^2$. As far as we know, the  B\"or\"oczky configuration as defined in \Cref{sec: borocky construction} is the only counterexample that does not contain points with multiplicity higher than three \cite[section 6, property 4]{PalastiFerudi}. Note this example of Elliptified B\"or\"oczky is not equivalent to one coming from the construction in \Cref{sec: borocky construction} as no planar B\"or\"oczky contains $48$ triple intersection points; the equation $\frac{n(n-3)}{6}+1=48$ has no integer solutions.
\begin{conjecture}
Let $6\mid n$ and consider any elliptic curve in $\PP^2$. Apply the Elliptified B\"or\"oczky construction to the points of $E[n]$. The defining ideal of $\mathcal{B}_{En}$, the triple points of the configuration, fails the containment $I^{(3)} \subseteq I^2$.
\end{conjecture}

\section{Open Problems}\label{sec: open problems}
We end this paper with a number of open problems, the first of which was the starting motivation of our investigation.

\begin{enumerate}
    \item Give a completely non-computer assisted proof of $I^{(3)}_n \not \subseteq I_n^2$.
    \item Can the B\"or\"oczky configuration, or the Elliptified B\"or\"oczky, be realized as a hyperplane arrangement with ordinary lines only? We could not do this even for \Cref{Fig: B12}. A positive answer would allow one to use the containment criteria developed in \cite{AlexandraBen}. 
    \item Investigate the construction in \Cref{sec: elliptified Boroczky} for larger values of $n$.
    \item Does applying the B\"or\"oczky construction to other elliptic curves, in higher dimensional projective spaces, also result in examples of ideals which fail the containment problem? 
    \item Give a non-computer assisted justification for the data in \Cref{tab: min gen data}.
\end{enumerate}
\appendix
\section{Macaulay2 Code}\label{section: M2 code}
\begin{small}
\begin{lstlisting}[language=Macaulay2]
K=toField(QQ[s,t]/ideal(t^2+t+1,s^3+2))
R=K[x,y,z]
--E6:= The 36 points that form the group Z_6 x Z_6 

E6={{{-1,1,0},{s,1,t},{-1,0,t},{t^2,t^2,s},{0,1,-t},{1,s,t}},
    {{1,s,1},{t^2,1,s},{s,t^2,t},{t,s,t},{t,t^2,s},{s,1,t^2}},
    {{0,1,-1},{t,s,1},{-t,1,0},{s,t^2,t^2},{-1,0,t^2},{t,1,s}},
    {{1,1,s},{s,t^2,1},{t,s,t^2},{t,t,s},{s,t,t^2},{t^2,s,1}},
    {{-1,0,1},{1,t,s},{0,1,-t^2},{t^2,s,t^2},{-t^2,1,0},{s,t,1}},
    {{s,1,1},{1,s,t^2},{t^2,t,s},{s,t,t},{t^2,s,t},{1,t^2,s}}};

--line:= Creates the lines of the Boroczky configuration.
--Lines:= Contains all configuration lines with duplicates.

line=(i,j)->
if det(matrix{{x,y,z},E6_i_j,E6_((3-2*i)%6)_((3-2*j)%6)})!=0 
then return det(matrix{{x,y,z},E6_i_j,E6_((3-2*i)%6)_((3-2*j)%6)}) 
else return x*(E6_i_j_0)^2+y*(E6_i_j_1)^2+z*(E6_i_j_2)^2;
Lines= 
flatten(for i from 0 to 5 list( for j from 0 to 5 list line(i,j))); 

--Defining auxiliary functions that remove repeated lines.
--dualPoint:= Reads the coefficients of a linear function.
--lineMeet:= Computes the intersection point of two lines.

dualPoint=i->{coefficient(x,i),coefficient(y,i),coefficient(z,i)};  
lineMeet=(i,j)-> 
dualPoint(det(matrix{{x,y,z},dualPoint(i),dualPoint(j)})); 
ND={};
for i from 0 to 34 do 
ND=append(ND,(select(Lines,j->lineMeet(j,Lines_i)=={0,0,0}))_0) 
Lines=unique(ND)

--Doubles:= Contains all points of multiplicity at least 2.
--Doubles contains (0,0,0) as well as duplicate points.
--intrank:= Counts how many lines meet at a point.

Doubles={};
for i from 0 to 17 do for j from i to 17 do
Doubles=append(Doubles,lineMeet(Lines_i,Lines_j));
Doubles=drop(unique(Doubles),1); 
NDpoints={};
for i from 0 to 137 do NDpoints=
append(NDpoints,(select(Doubles,j->rank(matrix{j,Doubles_i})==1))_0);
Doubles=unique(NDpoints); 
intrank=i->#select(Lines,j->sub(j,{x=>i_0,y=>i_1,z=>i_2})==0);  
Triples=select(Doubles,i->intrank(i)>2); 

--Checking that there are no Quadruple points.
--I is I_E6 defined in Section 5.

Quadruple=select(Doubles,i->intrank(i)>3);
I=intersect(apply(Triples,i->minors(2,matrix{{x,y,z},{i_0,i_1,i_2}})));
isSubset(saturate(I^3),I^2)
isSubset(ideal(product(Lines)),I^2)
   \end{lstlisting}   
\end{small}
\bibliographystyle{amsalpha}
\bibliography{biblio}
\end{document}